\documentclass[11pt]{amsart}
\usepackage{amsfonts}
\usepackage{amsmath,amsthm}
\usepackage{hyperref}
\usepackage{latexsym}
\usepackage{array}
\usepackage{amssymb}
\usepackage{enumerate}
\usepackage[notcite,notref]{showkeys}
\usepackage[francais,english]{babel}

\usepackage{color}
\usepackage[latin1]{inputenc}

\theoremstyle{plain}
\newtheorem{theorem}{Theorem}[section]
\newtheorem*{theorem*}{Theorem}
\newtheorem{lemma}[theorem]{Lemma}

\newtheorem{definition}[theorem]{Definition}

\theoremstyle{remark}
\newtheorem{remark}[theorem]{Remark}

\newtheorem*{lem*}{Lemma}
\newtheorem*{sublem*}{Sublemma}
\newtheorem*{remark*}{Remark}
\newtheorem*{NB*}{NB}

\newcommand{\R}{ \mathbb{R} }
\newcommand{\C}{ \mathbb{C} }
\newcommand{\Z}{ \mathbb{Z} }
\newcommand{\N}{ \mathbb{N} }

\newcommand{\T}{ \mathbb{T} }
\newcommand{\A}{ \mathbb{A} }

\newcommand{\fJ}{ \mathfrak J }

\newcommand{\an}{ \, \angle\,}

\newcommand{\cA}{ \mathcal{A} }
\newcommand{\cB}{ \mathcal{B} }
\newcommand{\cC}{ \mathcal{C} }

\newcommand{\Cc}{ \mathcal{C} }
\newcommand{\D}{ \mathcal{D} }

\newcommand{\E}{ \mathcal{E} }
\newcommand{\F}{ \mathcal{F} }

\newcommand{\J}{ \mathcal{J} }
\renewcommand{\L}{ \mathcal{L} }
\newcommand{\M}{ \mathcal{M} }
\newcommand{\cM}{ \mathcal{M} }

\newcommand{\NF}{ \mathcal{NF} }
\renewcommand{\O}{ \mathcal{O} }

\newcommand{\Tc}{ \mathcal{T} }
\newcommand{\cT}{ \mathcal{T} }

\newcommand{\om}{ \omega }
\newcommand{\Om}{ \Omega }
\newcommand{\ga}{\gamma }
\newcommand{\s}{ \sigma }
\newcommand{\De}{ (-\Delta) }
\newcommand{\ka}{ \kappa }
\renewcommand{\r}{ \rho }

\renewcommand{\phi}{ \varphi }
\newcommand{\eps}{\varepsilon}
\newcommand{\vark}{ \varkappa }
\newcommand{\de}{ \delta }

\newcommand{\zz}{\mathfrak z}
\renewcommand{\a}{ \alpha }

\newcommand{\LL}{ \Lambda^{\#}}

\newcommand{\la}{ \lambda_a }
\newcommand{\lb}{ \lambda_b }
\newcommand{\li}{ \lambda_i }
\newcommand{\lj}{ \lambda_j }
\newcommand{\lk}{ \lambda_k }
\newcommand{\lel}{ \lambda_\ell }
\newcommand{\yy}{ \rho }
\newcommand{\bH}{\mathbf{H} }
\newcommand{\La}{ \Lambda }

\newcommand{\diag}{\operatorname{diag}}
\newcommand{\Leb}{\operatorname{Leb}}

\newcommand{\meas}{\operatorname{meas}}

\newcommand{\lsim}{ \lesssim }

\def\ab#1{\left|#1\right|}
\def\aa#1{\left\Vert#1\right\Vert}

\newcommand{\be}{\begin{equation}}
\newcommand{\ee}{\end{equation}}
\newcommand{\ben}{\begin{equation*}}
\newcommand{\een}{\end{equation*}}
\newcommand{\ban}{\begin{align*}}
\newcommand{\ean}{\end{align*}}
\numberwithin{equation}{section}

\newcommand{\dd}{ \text{d} }

\newcommand{\p}{ \partial}

\author{L. Hakan Eliasson}
\address{Univ. Paris Diderot, Sorbonne Paris Cit\'e\\
Institut de Math\'emathiques de Jussieu-Paris rive gauche, UMR 7586\\
CNRS\\
Sorbonne Universit\'es, UPMC Univ. Paris 06\\
F-75013, Paris, France}
\email{hakan.eliasson@imj-prg.fr}

 \author{ Beno\^it Gr\'ebert}
\address{Laboratoire de Math\'ematiques Jean Leray, Universit\'e de Nantes, UMR CNRS 6629\\
2, rue de la Houssini\`ere \\
44322 Nantes Cedex 03, France}
\email{benoit.grebert@univ-nantes.fr}

\author{ Serge\"i B. Kuksin }
\address{CNRS\\
Institut de Math\'emathiques de Jussieu-Paris rive gauche, UMR 7586\\
Univ. Paris Diderot, Sorbonne Paris Cit\'e\\
Sorbonne Universit\'es, UPMC Univ. Paris 06\\
F-75013, Paris, France}
\email{sergei.kuksin@imj-prg.fr}

\title[KAM for multi-d PDE]
{A KAM theorem ~for space-multidimensional hamiltonian
 PDE}

\begin{document}

\begin{abstract}
We present an abstract KAM theorem, adapted to space-multidimensional hamiltonian
 PDEs with smoothing non-linearities. The main novelties of this theorem ~are that:
 \begin{itemize}
 \item the integrable part of the hamiltonian may contain a hyperbolic part and as a consequence the constructed invariant tori may be unstable.
 \item It applies to singular perturbation problem.
 \end{itemize}
 In this paper we state the KAM-theorem and
 comment on it, give the main ingredients of the proof, and present three applications of the theorem.
 \end{abstract}

\subjclass{ }
\thanks{
}
\maketitle

\tableofcontents

\section{Introduction}
In this paper we present and comment on an abstract KAM theorem,
 proved in \cite{EGK, EGK1}. In \cite{EGK} we focus on the main application of the theorem ~to
 the existence of small amplitude solutions for the nonlinear beam equation on a torus
 of any dimension, which was the
 motivation for establishing the theorem. In this short presentation we
 focus on the novelties of our result, give some elements of its proof and present
 two more applications of the theorem, in addition to that, described in \cite{EGK}.

\subsection{Notation}

We consider a Hamiltonian $H=h+ f$, where $h$ is a quadratic Hamiltonian
\be\label{h}
h=\Omega(\yy)\cdot r+\frac12 \sum_{a\in\L_\infty} \Lambda_a(\yy)(p_a^2+ q_a^2) +
\frac12\, \langle \bH(\yy)\zeta_\F,\zeta_\F\rangle\,.
\ee
Here
\begin{itemize}
\item $\rho$ is a parameter in $\D$, which is an open ball
 in the space $\R^{n}$;
\item $r\in\R^n$ are the actions corresponding to the internal modes
 $(r,\theta)\in (\R^n\times\T^n,dr\wedge d\theta)$;
\item $\L_\infty$ and $\F$ are respectively infinite and finite sets in $\Z^d$, $\L$ is the disjoint union $\L_\infty\cup\F$;
\item $\zeta=(\zeta_s)_{s\in\L}$ are the external modes, where
 $\zeta_s=(p_s,q_s)\in (\R^2, dq\wedge dp)$. The external modes decomposes in an infinite part $\zeta_\infty=
 (\zeta_s)_{s\in\L_\infty}$, corresponding to elliptic directions, and a finite part $\zeta_\F=(\zeta_s)_{s\in\F}$ which may contain hyperbolic directions;
\item the mappings
 \be\label{properties}\left\{\begin{array}{ll}
\Om:\D\to\R^{n},&\\
\La_a:\D\to \R\,,&\quad a\in \L_\infty,\\
\bH:\D\to gl({\F}\times {\F}),&\quad {}^t\! H=H,
\end{array}\right.\ee
are $\cC^{{s_*}}$-smooth, ${s_*}\ge 1$.
\item $f=f(r,\theta, \zeta;\rho)$ is the perturbation,
 small compare to the integrable part $h$ in a way, specified below in \eqref{epsest}.
\end{itemize}
The integrable Hamiltonian $h(r, \theta, \zeta_\F, \zeta_\infty)$ has a {\it finite-dimensional invariant torus}
\be\label{torus}
\{0\}\times\T^{n}\times\{0\}\times \{0\}\,,
\ee
and the equation, linearised on this torus, does not depend on
the angles $\theta$. This linearized equation
has infinitely many elliptic directions with purely imaginary eigenvalues
$$
\{\pm {\mathbf i}\Lambda(\r) : a\in \L_\infty\}
$$
and finitely many other directions given by the system
$$
\dot \zeta_\F= J\bH(\r)\zeta_\F
$$
(some of them may be hyperbolic).

\subsection{A perturbation problem}
The question which we address is if for most values of the parameter $\rho\in\D$
 the invariant torus \eqref{torus} persists under perturbations $h+f$ of
the Hamiltonian $h$, and, if so, if the perturbed equation, linearised
about solutions on this torus, is reducible to constant coefficients.

In finite dimension the answer to this question is affirmative
 under rather general conditions. For the first proof in the purely elliptic case see \cite{E88}, and for a more general case see \cite{Y99}. These statements say that, under general conditions, the invariant torus persists and remains reducible under sufficiently small perturbations, for a subset of parameters $\r$ of large Lebesgue measure. Since the unperturbed problem is linear, parameter selection can not be avoided here.

In infinite dimension the situation is more delicate, and results can only be proven under quite severe restrictions
on the set of
normal frequencies $\{\La_a\}$. In one space dimension these
 restrictions
are fulfilled for many PDEs; the first such result was obtained in \cite{Kuk87}.
For PDEs in higher space dimension the behaviour of the normal frequencies is much more complicated, and the
available results
are more sparse (see below).

Comparing with the existing results, the main novelties of the KAM theorem, stated in the next section are:
 \begin{itemize}
 \item It applies to singular perturbation problem, i.e. the size of the perturbation is coupled to the control that we have on the small divisors of the unperturbed part (see Subsection~\ref{sing}).
 \item The integrable part of the hamiltonian may contain a finite-dimensio\-nal hyperbolic part whose treatment requires higher smoothness in the parameters. If the hyperbolic part of the unperturbed linear system is non-trivial,
 the constructed invariant tori are unstable.
 \item
We have imposed no ``conservation of momentum'' on the perturbation. This allows to treat perturbations,
depending on the space-variable $x$, and
has the effect that during the KAM-iterations
our normal form is not diagonal in the purely elliptic directions.
In this respect it resembles the normal form, used in \cite{EK10} to treat the non-linear
Schr\"odinger equation.
\item A technical difference with previous works on KAM for PDE (including
 \cite{EK10}) is that now we use a different matrix norm with
 better multiplicative properties. This simplifies the functional analysis, involved in the proof.
\item
Comparing to \cite{EK10}, we impose a further decay property on the hessian of the non-quadratic part of the
Hamiltonian (see \eqref{reg}). As a consequence we do not
 have to use the involved
 T\"oplitz--Lipschitz machinery of the work \cite{EK10}. This simplifies the proof,
 but does not allow to apply the KAM theorem ~of this work
 to the NLS equations, unless we regularise the non-linearity as in Subsection~\ref{nls}.
 \end{itemize}

\subsection{Short review of the related literature}

If the KAM theory for 1d Hamiltonian PDEs is now well documented
 (see \cite{Kuk87,Kuk93,Pos89,Kuk00} for a first overview), only few results exist for the multidimensional equations.\\
Existence of quasi-periodic solutions of space-multidimensional PDE were first proved in \cite{B1} (see also
 \cite{B2}), using
 the Nash--Moser technic, which does not allow to analyse the linear stability of the obtained solutions.
 KAM-theorems which apply to some
 for small-amplitude solutions of multidimensional beam equations (see \eqref{beam} above) were obtained in \cite{GY1, GY2}. Both works treat equations with a constant-coefficient nonlinearity
 $g(x,u)=g(u)$, which is significantly easier than the general case.
 The first complete KAM theorem ~for space-multidimensional PDE was obtained in \cite{EK10}. Also see \cite{Ber1, Ber2}. \\
 The technic of the work \cite{EK10} has been developed
 in \cite{EGK, EGK1} to allow a KAM result without external parameters. There we proved the existence of
 small amplitude quasi-periodic solutions of the beam equation on the $d$-dimensional torus, investigated
 the stability of these solutions, and gave explicit examples of linearly unstable solutions, when the
 linearised equations have finitely many hyperbolic directions. These results are discussed in Section~3. \\
 NLS equations on the $d$-dimensional torus without external parameters were considered in \cite{WM} and \cite{PP1, PP2}, using the
 KAM-techniques of \cite{B1, B2} and \cite{EK10}, respectively. Their main
 disadvantage compare to the 1d theory (see \cite{KP}) is severe restrictions (in the form of
 a non-degeneracy condition) on the finite set of linear modes on which the quasi-periodic solutions are based. The
 notion of the non-degeneracy is not explicit so that it is not easy to give examples of
 non-degenerate sets of modes. \\
All these examples concern PDEs on the tori essentially because in that case the corresponding linear PDEs are
diagonalizable in the exponential basis and have rather specific and
similar spectral clusters.
 Recently some examples that do not fit this Fourier context have been considered: the Klein-Gordon equation on the sphere $\mathbb S^d$ (see \cite{GP1}) and the quantum harmonic oscillator on $\R^d$ (see \cite{GT} and \cite{GP2}). For the existence of quasi-periodic solutions for NLW and NLS on compact Lie groups via the
 Nash--Moser approach see \cite{BCP} and references quoted therein.

\

\section{Setting and statement of our KAM theorem}
In this section we state our KAM result for the Hamiltonian $H=h+ f$ as in the introduction.
\subsection{Setting} First of all we detail the structures behind the objects appearing in \eqref{h} and the hypothesis needed for the KAM result.

\noindent {\bf Linear space.}
For any $\ga=(\ga_1,\ga_2)\in\R^2$ we denote by $Y_\gamma$ the following weighted { complex}
 $\ell_2$-space
\be\label{Y}
Y_\ga=
\Big\{\zeta= { \Big(
\zeta_s = \left(\begin{array}{ll}\xi_s\\ \eta_s \\ \end{array}\right)
\in\C^2,}\ s\in \L\Big) \mid \|\zeta\|_\ga<\infty\Big\} ,
\ee
where\footnote{We recall that $|\cdot|$ signifies the Euclidean norm.}
$$
\|\zeta\|_\ga^2=\sum_{s\in\L}|\zeta_s|^2\langle s\rangle^{2\ga_2}e^{2\gamma_1 |s|},\qquad
\langle s\rangle= \max (|s|,1).
$$
Endowed with this norm, $Y_\ga$ is a Banach space. Furthermore if $\ga_2>d/2$, then this space is an algebra with respect to the convolution.

 In a space $Y_\gamma$ we define the complex conjugation as the involution
\be\label{inv}
\zeta={}^t(\xi, \eta)\mapsto {}^t(\bar\eta, \bar\xi)\,.
\ee
Accordingly, the real subspace of $Y_\gamma$ is the space
\be\label{reality}
Y_\ga^R=
\Big\{
\zeta_s = \left(\begin{array}{ll}\xi_s\\ \eta_s \\ \end{array}\right)\mid \eta_s=\bar\xi_s,
\ s\in \L\Big\} \,.
\ee
Any mapping defined on (some part of) $Y_{\ga}$ with values in a complex Banach space with
a given real part
is called {\it real} if it gives real values to real arguments.

\medskip

 \noindent {\bf Infinite matrices.}
Let us define the pseudo-metric on $\Z^{d_*}$
$$
(a,b)\mapsto [a-b]=\min (|a-b|,|a+b|).$$
We shall consider matrices $A:\L\times\L\to gl(2,\C)$, formed by $2\times2$-blocs
(each $A_a^{b}$ is a $2\times2$-matrix). Define
\be\label{matrixnorm}
|A|_{\ga,\vark}=\max\left\{\begin{array}{l}
\sup_a\sum_{b} \ab{A_a^b} e_{\ga,\vark}(a,b)\\
\sup_b\sum_{a} \ab{A_a^b} e_{\ga,\vark}(a,b),
\end{array}\right.\ee
where the norm on $A_a^{b}$ is the matrix operator norm and where the weight $e_{\ga,\vark}$ is defined
by
\be\label{weight}
e_{\ga,\vark}(a,b)=Ce^{\ga_1[a-b]}\max([a-b],1)^{\ga_2}\min(\langle a\rangle,\langle b\rangle)^\vark\ee
for any\footnote{\ $(\ga_1',\ga_2')\le (\ga_1,\ga_2)$ if, and only if $\ga_1'\le\ga_1$ and $\ga_2'\le \ga_2$} $\ga=(\ga_1,\ga_2)\ge(0,0)$, $\vark\ge0$ and for some constant $C$ depending on $\ga$, $\vark$.

Let $\cM_{\ga,\vark}$ denote the space of all matrices $A$ such that
$\ab{A}_{\ga,\vark}<\infty$. Clearly $\ab{\cdot}_{\ga,\vark}$ is a norm
on $\cM_{\ga,\vark}$. It follows by well-known results that $\cM_{\ga,\vark}$,
provided with this norm, is a Banach space. Compare to the $\ell^\infty$-norm used in \cite{EK10}, the $\ell^1$-norm \eqref{matrixnorm} has the great advantage to enjoye, when $\ga_2\geq \vark$, the algebra property
$$\ab{{BA} }_{\ga,\vark}\le \ab{A}_{\ga,0} \ab{B}_{\ga,\vark}$$
and to satisfy
$$
\aa{A\zeta }_{\tilde \ga}\le \ab{A}_{\ga,\vark}\aa{\zeta }_{\tilde\ga}\,,$$
if $-\ga\le \tilde\ga\le\ga$.
In particular, for any $-\ga\le \tilde\ga\le\ga$, we have a continuous embedding of $\cM_{\ga,\vark}$,
$$\cM_{\ga,\vark}\hookrightarrow \cM_{\ga,0}\to
 \cB(Y_{\tilde \ga},Y_{\tilde \ga}),$$
into the space of bounded linear operators on $Y_{\tilde \ga}$. Matrix multiplication in $\cM_{\ga,\vark}$
corresponds to composition of operators.

\smallskip

For our applications we must consider a larger sub algebra with somewhat weaker decay properties.
For $\ga=(\ga_1,\ga_2)\ge (0,m_*)$ with $m_*>d/2$ fix,
let
\be\label{b-space}
\cM_{\ga,\vark}^b= \cB(Y_{\ga},Y_{\ga})\cap \cM_{(\ga_1,\ga_2-m_*),\vark}\ee
which we provide with the norm
\be\label{b-matrixnorm}
\aa{A}_{\ga,\vark }=\aa{A}_{ \cB(Y_{\ga},Y_{\ga})}+ \ab{A}_{(\ga_1,\ga_2-m_*),\vark }.\ee
This norm makes $\cM_{\ga,0}^b$ into a Banach sub-algebra of $\cB(Y_{\ga};Y_{\ga})$
and $\cM_{\ga,\vark}^b$ becomes an ideal in $\cM_{\ga,0}^b$.

\medskip

\noindent {\bf A class of Hamiltonian functions.}
Let
$$\ga=(\ga_1,\ga_2)\ge (0,m_*+\vark)=:\ga^*.$$
For a Banach space $B$ (real or complex) we denote
$$
\O_s(B)=\{x\in B\mid \|x\|_B<s\}\,,
$$
and for $\sigma,\ga,\mu\in(0,1]$ we set
\begin{align*}
\T^n_\s=&\{\theta\in\C^n/2\pi\Z^n\mid |\Im \theta|<\sigma\},\\
\O^\ga(\s,\mu)=& \O_{\mu^2}(\C^n) \times \T^n_\s \times \O_\mu(Y_\ga)=\{(r, \theta, \zeta)\}\,.
\end{align*}
We will denote the points in $\O^\ga(\s,\mu)$ as $x=(r,\theta,\zeta)$. \\

 Fix $s^*\geq 0$. Let $f:\O^{\ga^*}(\s,\mu)\times\D\to \C$ be a $C^{s^*}$-function, real holomorphic
 in the first variable $x=(r,\theta,\zeta)$,
 such that for all $0\le\ga'\le\ga$ and all $\yy\in\D$ the gradient-map
 $$
 \O^{\ga'}(\s,\mu)\ni x\mapsto \nabla_\zeta f(x,\yy)\in Y_{\ga'}$$
 and the hessian-map
 $$
 \O^{\ga'}(\s,\mu)\ni x\mapsto \nabla^2_{\zeta} f(x,\yy)\in \M^b_{\ga',\vark}
 $$
 also are real holomorphic.
 We denote this set of functions by ${\Tc}_{\ga,\vark,\D}(\s,\mu)$.
 For a function $h\in {\Tc}_{\ga,\vark,\D}(\s,\mu)$ we define
 the norm
 $$|h|_{\substack{\s,\mu\ \ \\\ga,\vark,\D}}$$
 through
\be\label{schtuk}
\sup_{
\substack{
0\le\ga'\le \ga\\ j=0,\cdots,s^*
}}\
\sup_{
\substack{
x\in O^{\ga'}(\s,\mu)\\ \yy\in\D
}}
\max( |\partial^j_\yy h(x,\yy)|,
{ \mu} \|\partial^j_\yy \nabla_\zeta h(x,\yy)\|_{\ga'},
{\mu^2}\|\partial^j_\yy \nabla^2_\zeta h(x,\yy)\|_{\ga',\vark}).
\ee

 \begin{remark}\label{decay} We note that if $\vark>0$, then even the diagonal of the hessian of
 $f\in {\Tc}_{\ga,\vark,\D}(\s,\mu)$
 has a decay property since then
 \be\label{reg} |\nabla^2_{\zeta_a,\zeta_b} f|\leq C\frac{e^{-\ga_1[a-b] }}{\langle a\rangle^{\vark}\langle b\rangle^{\vark}}.\ee
 This will be crucial to preserve the second Melnikov property (see Assumption~A3 above and Subsection~\ref{mel}) during the KAM iterations.
 \end{remark}

 For any function $h\in {\Tc}_{\ga,\vark,\D}(\s,\mu)$ we denote by $h^T$ its Taylor polynomial
at $r=0, \zeta=0$, linear in $r$ and quadratic in $\zeta$:
$$
h(x,\yy)=h^T(x,\yy)+ O(|r|^2+\|\zeta\|^3+|r|\|\zeta\|).
$$

\subsection{KAM Theorem}
We consider the Hamiltonian $H=h+f$ with $f\in {\Tc}_{\ga,\vark,\D}(\s,\mu)$ and $h$ as in \eqref{h},
and assume that $h$ satisfies the Assumptions~A1 -- A3, depending on constants
\be\label{const}
 \delta_0, c, C, \beta> 0,\;\; s_*\in\N\,.
\ee
 To formulate the assumptions we first introduce the partition of $ \L_{\infty}$ to the clusters
 $[a]$, given by
\be\label{partition1}[a]=\left\{\begin{array}{ll}
\{b\in\L_\infty: |b|\le c\} & \ \textrm{if}\ \ab{a}\le c\\
\{b\in\L_\infty: |b|=|a|\}& \ \textrm{if}\ \ab{a} > c,\end{array}\right.\ee
where $c$ is some (possibly quite large) constant.

\smallskip

\noindent{\bf Hypothesis A1} (spectral asymptotic.)
 For all $\yy\in\D$ we have \\

(a)\; $|\Lambda_a|\ge \delta_0$ $\;\ \forall\, a\in \L_\infty$;\\

(b) \; $| \Lambda_a-|a|^{2}|\le c \langle a\rangle^{-\beta}$ $\;\ \forall\, a\in \L_\infty$;\\

(c) \; $\|(J\bH (\yy))^{-1}\| \le \frac1{\delta_0}\,,\;\;\ \|(\Lambda_a(\yy) I -iJ \bH(\yy))^{-1}\| \le \frac1{\delta_0}\;\;\ \forall\, a\in \L_\infty\,;
$\\

(d)
 $|\Lambda_a(\rho) +\Lambda_b(\rho)| \ge \delta_0$ for all $a,b\in \L_\infty$;
 \\

(e)
 $|\Lambda_a(\rho) -\Lambda_b(\rho)| \ge \delta_0$ if $a,b\in \L_\infty$ and $[a]\ne[b]$.
\bigskip

\noindent
{\bf Hypothesis A2} (transversality).
 For each $k\in \Z^n\setminus \{0\}$ and every vector-function
$\Omega'(\yy)$ such that $|\Omega'-\Omega|_{C^{s_*}(\D)}\le \delta_0$ there exists
a unit vector $\zz=\zz(k)\in\R^n$, satisfying
\be\label{o}
|\p_\zz \langle k, \Omega'(\rho)\rangle |\ge \delta_0\qquad \forall\, \rho\in\D.
\ee
Besides the following properties (i)-(iii) hold for each $k\in \Z^n\setminus \{0\}$:

(i) For any $a,b \in\L_\infty\cup \{\emptyset\}$ such that $(a,b)\ne (\emptyset, \emptyset)$,
consider the following operator, acting on the space of $[a]\times[b]$-matrices \footnote{so if $b=\emptyset$, this is the space
$\C^{[a]}$.}
$$
L(\yy): X\mapsto ({\Omega'}(\yy)\cdot k)X \pm Q(\yy)_{[a]} X
+ X Q(\rho)_{[b]}\,.
$$
Here $Q(\yy)_{[a]}$ is the diagonal matrix diag$\{\Lambda_{a'}(\rho) : a'\in[a]\}$,
 and $Q(\yy)_{[\emptyset]}=0$. Then either
\be\label{invert}
\|L(\rho)^{-1}\|\le \delta_0^{-1}\qquad \forall\, \rho\in\D\,,
\ee
or there exists a unit vector $\mathfrak z$ such that
$$
|\langle v,\p_\zz L(\rho)v\rangle| \ge\delta_0\qquad \forall\, \rho\in\D\,,
$$
for each vector $v$ of unit length.

(ii) Denote $m=2|\F|$ and consider the following operator in $\C^m$, interpreted as a space of
row-vectors:
$$
L(\yy,\lambda): X\mapsto ({\Omega'}(\yy)\cdot k)X +\lambda X+i XJ\bH(\yy)\,.
$$
 Then
$$
\|L^{-1}(\yy,\Lambda_a(\r))\| \le\delta_0^{-1}\qquad \forall\, \yy\in\D,\;\; a\in\L_\infty\,.
$$

(iii) For any $a,b \in \F \cup \{\emptyset\}$ such that $(a,b)\ne (\emptyset, \emptyset)$,
consider the operator, acting on the space of $[a]\times[b]$-matrices:
$$
L(\yy): X\mapsto ({k\cdot\Omega'}( \yy))X - iJ\bH(\yy)_{[a]}X + iXJ\bH(\yy)_{[b]}
$$
(the operator $\bH(\yy)_{[a]}$ equals $\bH(\rho)$ if $a\in\F$ and equals 0 if $a=\emptyset$, and similar
with $\bH(\yy)_{[b]}$).
Then the following alternative holds: either $L(\yy)$ satisfies \eqref{invert}, or there exists an integer $1\le j\le s_*$
such that
\be\label{altern1}
|\p_\zz ^j\det L(\rho)| \ge \delta_0 \|L(\rho)\|_{C^j(\D)}^{m -2}\qquad \forall\, \rho\in \D\,.
\ee
Here $m=4|\F|^2$
 if $a,b\in\F$ and $m=2|\F|$ if $a$ or $b$ is the empty set.

\bigskip

\noindent
{\bf Hypothesis A3} (the Melnikov condition). There exist $\tau>0$, ${\yy_*}\in\D$
and $C>0$ 
such that
\be\label{melnikov}
| k\cdot\Om(\yy_*) -(\Lambda_a({\yy_*}) - \Lambda_b({\yy_*}))| \ge C|k|^{-\tau}\quad
\forall\, k\in\Z^n, k\ne0,\;\text{if}\;
 a,b\in \L_\infty \setminus[0].
 \ee

 Denote
$$
\chi = | \p_\yy\Omega(\yy)|_{C^{s_*-1}}
+\sup_{a\in\L_\infty} |\p_\yy\La(\yy)|_{C^{s_*-1}} +
\|\p_\yy \bH\|_{C^{s_*-1}} \,.
$$
Consider the perturbation $f(r,\theta, \zeta; \yy)$ and assume that
$$
\eps=|f^T|_{\substack{\s,\mu\ \ \\\ga,\vark,\D}}<\infty\,,\quad \xi =|f|_{\substack{\s,\mu\ \ \\\ga,\vark,\D}}<\infty\,,
$$
for some $\gamma, \sigma, \mu\in(0,1]$.
 We are now in position to state the abstract KAM theorem ~from \cite{EGK}.
 More precisely, the result below follows from
 Corollary~6.9 of \cite{EGK}.

\begin{theorem}\label{main} Assume that Hypotheses~A1-A3 hold for $\yy\in\D$. Then there exist
$ \a, C_1>0$ and $\eps_*=\eps_*(\chi,\xi,\delta_0)>0$
such that if
\be\label{epsest2}
\eps \le \eps_*(\chi,\xi,\delta_0)\,,
\ee
then
there is a Borel set
$\D'\subset \D$ with
$$\
\text{meas}(\D\setminus \D')\leq C_1\eps^{\a}$$
and there exists a $C^{s_*}$-smooth mapping
 $$
 {\mathfrak F} :\O^{\ga^*}(\s/2,\mu/2)\times\D
 \to \O^{\ga^*}(\s,\mu)\,,\quad (r,\theta,\tilde\zeta;\yy)\mapsto {\mathfrak F}_\yy(r,\theta,\tilde\zeta)\,,
 $$
 defining for $\rho\in\D$
 real holomorphic symplectomorphisms
${\mathfrak F}_\yy :\O^{\ga^*}(\s/2,\mu/2)\to \O^{\ga^*}(\s,\mu)$, satisfying for any $x\in \O^{\ga^*}(\s/2,\mu/2)$, $\rho\in \D$ and $|j|\le s_*$
the estimates
\be\label{Hest}
\| \p_\yy^j ({\mathfrak F}_\yy(x) -x) \|_{0} \le C_1 {\frac{\eps}{\eps_*}}\,, \qquad
\| \p_\yy^j ( d{\mathfrak F}_\yy (x)- I) \|_{0,0} \le C_1 {\frac{\eps}{\eps_*}}\,,
\ee
 such that for $\r\in\D$
\be\label{nf}
H\circ {\mathfrak F}_\yy=
\tilde \Om(\yy)\cdot r +\frac 1 2 \langle \zeta, A(\yy)\zeta\rangle +g(r,\theta,\zeta;\yy).
\ee
and for $\r\in\D'$
\be\label{invari}
\partial_\zeta g=\partial_r g=\partial^2_{\zeta\zeta}g=0\; \;\text{for}\;\;
\zeta=r=0\;.
\ee
 Here $\tilde\Om=\tilde\Om(\yy)$ is a new frequency vector satisfying
 \be\label{estimOM}
 \|\tilde\Om-\Om\|_{\mathcal C^{s_*}}\leq C_1 {\frac{\eps}{\eps_*}} \,,
 \ee
 and $ A:\L\times\L \to \M_{2\times 2}(\yy)$
 is an infinite real symmetric matrix belonging to
 $\M_{\ga^*,\vark}^b$. It is of the form $A=A_f\oplus A_\infty$, where
 \be\label{k6}
\| \p_\yy^\alpha(
 A_f(\yy) - \bH(\yy) )\|\le C_1 {\frac{\eps}{\eps_*}} \,,\quad |\alpha|\le s_*-1\,.
\ee
The operator $A_\infty$ is such that
 $(A_{\infty})_{ a}^b=0$ if $[a]\ne[b]$ (see \eqref{partition1}), and
 all eigenvalues of the hamiltonian operator $JA_\infty$ are pure
 imaginary.

\end{theorem}

So for $\yy\in\D'$ the torus ${\mathfrak F}_\yy \big(\{0\}\times\T^n\times\{0\}\big)$
is invariant for the hamiltonian system with
the Hamiltonian $H(\cdot;\yy)=h+f$ with $h$ given by \eqref{h}, and the hamiltonian
 flow on this torus is conjugated by the map ${\mathfrak F}_\yy$ with the
linear flow, defined by the Hamiltonian \eqref{nf} on the torus $(\{0\}\times\T^n\times\{0\})$.

\begin{remark}\label{rem-sing}
Estimate \eqref{epsest2} is the crucial assumption: it links the size of the perturbation with the Hypotheis A2-A3 on the unperturbed part. In particular $\eps_*$ depends on $\delta_0$ and since the perturbation has to be negligible compare with the control that we have on the small divisors, we can expect $$\eps \ll \delta_0\ .$$
We will see in Theorem \ref{main2} that the link is more involved and also depends on the size of $\chi$ and $\xi$.\end{remark}

%
%
%
%
%

\section{Elements of the proof}
The proof has the structure of a classical KAM-theorem ~carried out in a complex infinite-dimensional situation. The main part is,
as usual, the solution of the homological equation with reasonable estimates. The fact that the block structure is not
diagonal complicates, but this was also studied in for example \cite{EK10} (see Section~\ref{bloc}). The iteration combines a finite linear iteration
with a ``super-quadratic'' infinite iteration (see Section~\ref{iter}). This has become quite common in KAM and was also used
in \cite{EK10}.\\
In this section we also focus on the new ingredients: the use of the non linear homological equation that leads to better
estimates than the standard ones (see Section~\ref{hom}); the decay property \eqref{reg} which is very
useful to preserve the second Melnikov property during the KAM iteration (see Section~\ref{mel}); the treatment of the hyperbolic directions (see Section~\ref{hype}).

\subsection{Block decomposition, normal form matrices} \label{bloc}
In this subsection we recall two notions introduced in
\cite{EK10} for the nonlinear Schr\"odinger equation.
They are essential to overcome the problems of small divisors
in the multidimensional context.

 \medskip

\noindent{\bf Partitions.}
 For any $\Delta\in\N\cup \{\infty\}$
we define an equivalence relation on $\Z^{d_*}$, generated by the pre-equivalence relation
$$ a\sim b \Longleftrightarrow \left\{\begin{array}{l} |a|=|b| \\ {[a-b]}
 \leq \Delta. \end{array}\right.$$
Let $E_\Delta(a)$ denote the equivalence class of $a$ -- the {\it block} of $a$.
The crucial fact is that the blocks have a finite maximal diameter
$$d_\Delta=\max_{E_\Delta(a)=E_\Delta(b)} [a-b]$$
which do not depend on $a$ but only on $\Delta$:
\be\label{crucial}
d_\Delta\leq C \Delta^{\frac{(d_*+1)!}2}.
\ee
 This was proved in \cite{EK10}.

If $\Delta=\infty$ then the block of $a$ is the sphere $\{b: |b|=|a|\}$.
Each block decomposition is a sub-decomposition of the trivial decomposition formed by the spheres $\{|a|=\text{const}\}$.

\smallskip

On $\L_{\infty}\subset\Z^{d_*} $ we define the partition
$$[a]_\Delta=
\left\{\begin{array}{ll}
E_\Delta(a)\cap\L_{\infty} & \ \textrm{if}\ a\in \L_{\infty}\ \textrm{and}\ \ab{a}>C\\
 \{ b\in \L_{\infty}: \ab{b}\le C\} & \ \textrm{if}\ a\in \L_{\infty} \ \textrm{and}\ \ab{a}\le C\end{array}\right.$$
-- when $\Delta=\infty$, then this is the partition \eqref{partition1}. We extend this partition on $\L=\F\sqcup\L_{\infty}$ by setting $[a]_\Delta=\F$ if $a\in\F$. We denote it $\E_\Delta$.

\medskip

\noindent{\bf Normal form matrices.}
If $A:\ \L\times \L\to gl(2, \C)$ we define its {\it block components}
$$
A_{[a]}^{[b]}:[a]\times[b]\to gl(2, \C)$$
to be the restriction of $A$ to $[a]\times[b]$. $A$ is {\it block diagonal} over $\E_\Delta$ if, and only if,
$A_{[b]}^{[a]}=0$ if $[a]\neq [b]$. Then we simply write $A_{[a]}$ for $A_{[a]}^{[a]}$.

On the space of $2\times 2$ complex matrices we introduce a projection
$$
\Pi: gl(2, \C)\to \C I+\C J,
$$
 orthogonal with respect to the Hilbert-Schmidt scalar product. Note that
$\C I+\C J$ is the space of matrices, commuting with the symplectic matrix $J$.
\begin{definition}\label{d_31}
We say that a matrix $A:\ \L\times \L\to gl(2, \C)$ is on normal form with respect to
$\Delta$, $\Delta\in\N\cup \{\infty\}$, and write $A\in \NF_\Delta$, if
 \begin{itemize}
 \item[(i)] $A$ is real valued,
 \item[(ii)] $A$ is symmetric, i.e. $A_b^a\equiv {}^t\hspace{-0,1cm}A_a^b$,
 \item[(iii)] $A$ is block diagonal over $\E_\Delta$,
 \item[(iv)] $A$ satisfies $\Pi A^a_b\equiv A^a_b$ for all $a,b\in\L_{\infty}$.

 \end{itemize}
 \end{definition}

By extension we say that a Hamiltonian is on normal form if it reads
\be\label{hnf}
h(r,\zeta,\r)=\Om'(\yy)\cdot r +\frac 1 2 \langle \zeta, A(\yy)\zeta\rangle\ee
with $A$ a matrix in normal form and $\Om'$ close to $\Om$ in $C^1$ norm on $\D$.

The real quadratic form ${\mathbf q}(\zeta)= \frac 1 2\langle \zeta,A\zeta \rangle$, $\zeta=(p,q)$,
reads
$$
\frac 1 2\langle p,A_{1}p \rangle+\langle p,A_{2}q \rangle+\frac 1 2\langle q,A_{1}q \rangle
+\frac 1 2\langle \zeta_{\F}, H(\r) \zeta_{\F}\rangle
$$
where $A_{1}$ and $H$ are real symmetric matrices and $A_{2}$ is a real skew symmetric matrix. Note that in the complex variables $z_a=(\xi_a,\eta_a)$ defined through
$$
\xi_a=\frac 1 {\sqrt 2} (p_a+{\mathbf i}q_a),\quad \eta_a =\frac 1 {\sqrt 2} (p_a-{\mathbf i}q_a),$$
for $a\in\L_\infty$, and acting like the identity on $ (\C^2)^\F$,
the quadratic form ${\mathbf q}$ reads
$$
\langle \xi,Q\eta\rangle
+\frac 1 2\langle z_{\F}, H(\r) z_{\F}\rangle ,$$
where
$$Q= A_{1}+{\mathbf i}A_{2}.$$
Hence $Q$ is a Hermitian matrix.

The value of $\Delta$ will grow during the KAM iteration. At the beginning the Hamiltonian $h$ given in \eqref{h} is in normal form with respect to $\E_\Delta$ for any $\Delta\geq 1$. At the end of the story, i.e. in \eqref{nf}, $h_\infty=\tilde \Om(\yy)\cdot r +\frac 1 2 \langle \zeta, A(\yy)\zeta\rangle$ is in normal form with respect to $\E_\infty$.

\subsection{Homological equation}\label{hom}
Let us first recall the general KAM strategy.
Let $h$ be a the Hamiltonian given in \eqref{h}. Let $f$ be a perturbation and
 $$f^T=f_\theta+\langle f_r, r\rangle+\langle f_\zeta,\zeta\rangle+\frac 1 2 \langle f_{\zeta\zeta}\zeta,\zeta \rangle$$
 be its jet. The torus $\{0\}\times\T^n\times\{0\}$ is a KAM torus (i.e. an invariant torus on which the angles move linearly) for the Hamiltonian $h$. We want to prove the persistency of this KAM torus, in a deformed version, for $H=h+f$.
 If $f^T$ equals zero, then $\{0\}\times\T^n\times\{0\}$ would be still invariant by the flow generated by
$h+f$ and we were done. In general we only know that $f^T$ is small, say $f^T=\O (\eps)$. In order to
decrease the error term
 we search for a hamiltonian jet $S=S^T=\O (\eps)$ such that its time-one flow map
$\Phi_S=\Phi_S^1$ transforms the Hamiltonian $h+f$ to
$$
(h+f)\circ \Phi_S=h^+ + f^+,
$$
where $h^+$ is a new normal form, $\eps$-close to $h$, and the new perturbation $f^+$ is such that its jet
is much smaller than $f^T$. More precisely,
$$
h^+=h+\tilde h,\qquad
\tilde h=c(\r)+\langle \chi(\r),r\rangle+ \frac 1 2 \langle \zeta, B(\r)\zeta\rangle=\O(\eps),
$$
with $B$ on normal form and
$\
\left(f^+\right)^T=\O (\eps^2).
$

As a consequence of the Hamiltonian structure we have (at least formally) that
$$(h+f)\circ \Phi_S= h+\{ h,S \}+\{ f-f^T,S \}+f+ \O (\eps^2).$$
So to achieve the goal above
we should solve the {\it nonlinear homological equation}\footnote{The equation is nonlinear, because the solution $S$ depends nonlinearly on $f$.}:
\be \label{eq-homonl}
\{ h,S \}+ \{ f-f^T,S \}^T+f^T=\tilde h.
\ee
Then we repeat
the same procedure with $h^+$ instead of $h$ and $f^+$ instead of $f$. Thus we will have to solve the homological equation, not
only for the normal form Hamiltonian \eqref{h}, but for
more general normal form Hamiltonians
\eqref{hnf} with $\Om'$, $\eps-$close to $\Om$, and $A$ in normal form and $\eps-$close to $A_0=\diag(\La_a,\ a\in\L_\infty)\oplus\bH$.

\medskip

\noindent{\bf Nonlinear homological equation versus homological equation:} In many proofs of KAM theorems, one uses the homological equation
\be \label{eq-homo}
\{ h,S \}+f^T=0
\ee
instead of the nonlinear one \eqref{eq-homonl}. In that case we have an extra term in the jet of $f^+$ which is $\{ f-f^T,S \}^T$. At the first step of the iteration we have $f=\O (\eps)$ thus this is a term of order $\eps^2$ and this is not a problem. Nevertheless, although at each step $f^T$ is smaller, this is not the case for $f$ which remains of order $\eps$. So at step $k$ $\{f_k,S_k\}$ is of order $\eps \eps_k^2$ and not of order $\eps_k^2$. This problem can be overcome by using that $f(x)-f^T(x)=O(\|x\|^3)$ and thus is small for $\|x\|$ small. But this imposes a important constraint on the size $\mu_k$ of the analyticity domain of the Hamiltonian, $f_k\in {\Tc}_{\ga_k,\vark,\D_k}(\s_k,\mu_k)$. Essentially we have to choose ${\mu_k}\leq \eps_k^{\alpha}$ for some $\alpha>0$ (see for instance \cite{EK10}). The counterpart of such a choice is paid each time we use Cauchy's estimates in the variables $r$ or $\zeta$. In the present work this would drastically modified the condition \eqref{epsest} and as a consequence the theorem ~would not be sufficiently efficient to deal with a singular perturbation problem as the one presented in the subsection \ref{sing}.

\subsection{Small divisors and Melnikov condition}\label{mel}
The KAM proof is based on an iterative procedure that requires
to solve a homological equation at each step. Roughly speaking, it consists in inverting an infinite dimensional matrix whose eigenvalues are the so-called small divisors:
\begin{align*}
& \om\cdot k\quad k\in\Z^\A,\\
&\om\cdot k +\lambda_a\quad k\in\Z^\A,\ a\in\L,\\
&\om\cdot k +\lambda_a\pm\lambda_b\quad k\in\Z^\A,\ a,b \in\L
\end{align*}
where $\om=\om(\r)$ and $\la=\la(\r)$ are small perturbations (changing at each KAM step) of the original frequencies $\Om(\r)$ and $\Lambda_a(\r)\ a\in\L$ the eigenvalues of $A_0=\diag(\Lambda_a,\ a\in\L_\infty\})\oplus\bH$. In this subsection we focus on the elliptic part of $A_0$, i.e. on the case $a,b\in\L_\infty$.\\
Ideally we would like to bound away from zero all these small divisors.
In particular, this leads to infinitely many non resonances conditions of the type
$$|\om\cdot k +\la-\lb|>\frac{\ka}{|k|^\tau}, \quad k\in\Z^\A,\ a,b \in\L_\infty$$
for some parameters $\ka>0$ and $\tau>0$. Of course we have to exclude the case $k=0,\ a=b$ for which the small divisor is identically zero and this is precisely the reason why the external frequencies $\la$, $a\in\L_\infty$, move at each step.\\
When $d=1$ we have $|\la-\lb|\geq 2|a|$ for $|b|\neq|a|$. Therefore for each fixed $k$ there are only finitely many non resonances conditions and we can expect to satisfy them for a large set of parameters $\r$.\\
Now when $d\geq 2$, the frequencies $\la,\ a\in\L_\infty$, are not sufficiently separated and we really have to manage infinitely many non resonances conditions for each $k$. In general, it is not possible to control so many small divisors. Part of the solution consists in decomposing $\L$ in blocks $[a]_\Delta$ and to solve the homological equation according to this clustering. Then we only have to control the small divisors
$$|\om\cdot k +\la-\lb| \quad \text{ for }k\in\Z^\A,\ a,b \in\L_\infty,\ [a]_\Delta\neq[b]_\Delta $$
which is more reasonable. Actually when $|a|=|b|$ then $[a]_\Delta\neq[b]_\Delta$ implies $[a-b]\geq \Delta$. At this stage we have to recall that we want to control the small divisor $|\om\cdot k +\la-\lb|$ precisely to kill the quadratic term of the perturbation $\partial^2_{\xi_a\eta_b}f(\theta,0,0) \xi_a\eta_b$. But when $[a-b]\geq \Delta$, we can use the off diagonal exponential decay of the corresponding Hessian term in $\mathcal T_{\ga,\vark,\D}(\s,\mu)$ (see \eqref{matrixnorm} and \eqref{weight}) to assert that
\be\label{yes}\partial^2_{\xi_a\eta_b}f (\theta,0,0)=O(e^{-\ga\Delta})
\ee
 i.e. this term is already very small and it is not necessary to kill it.

Then it remains to consider the case where $|a|\neq |b|$. In that case, by hypothesis A3,
we can control from below $|\Om\cdot k +\La_a-\La_b|$ for all $k\neq 0,$ $a,b\in\L_\infty$. Then we get
$$|\om\cdot k +\La_a-\La_b|\geq \ka$$
 for all $a,b\in\L_\infty$ and $|k|$ not too large compared to $\ka^{-1}$ since $\om$ is close from $\Om$. On the other hand, as a consequence of the decay property \eqref{reg}, we can verify that
$$|\la(\rho)-\La_a(\rho)|\leq \frac{C}{|a|^{2\vark}}.$$
Therefore, the control that we have on $|\om\cdot k +\La_a-\La_b|$ leads to the control of $|\om\cdot k +\la-\lb|$ for $a$ and $b$ large enough (depending on $|k|$).\\ Now if $a$ or $b$ is small, says less that $M$, the other one has to be less than $C|k|+M$ in such a way $|\la-\lb|$ is comparable to $\om\cdot k$ and the small divisor can be small. At the end of the day, at fix $k$, it remains to control only finitely many small divisor and this can be achieved excising the possibly wrong subset of parameters.
By hypothesis A2(i), 
we always have a direction in which the derivative of the small divisor is larger than $\delta_0$ and thus the excised subsets are
small.

\smallskip

Finally we note that since the size of the block $[a]_\Delta$ does not depend
on the size of the index $a$ but only on $\Delta$ (see \eqref{crucial}) all the norms of $[a]_\Delta\times[b]_\Delta$ matrices are equivalent modulo constants that only depend on $\Delta$. Thus we can solve the homological equation in infinity matrix norm and then we can deduce estimates in operator norm.


\subsection{Small divisor and hyperbolic part}\label{hype}
Let us now consider small divisors involving the hyperbolic part. We will focus on the control of
\be\label{hy}\om\cdot k +\lambda_a\quad k\in\Z^\A,\ a\in\F.\ee
As in the previous section we want to solve the homological equation according to the clustering $\E_\Delta$. This means that, instead of trying to control the small divisors \eqref{hy} for each $a\in\F$, we want to control the inverse of the matrix
$$L(\r)=\om(\r)\cdot k I+iJ\bH'(\r)$$
where $\om=\om(\r)$ and $\bH'=\bH'(\r)$ are small perturbations (changing at each KAM step) of $\Om(\r)$ and $\bH(\r)$. The difference with the previous section is that now we are not dealing with Hermitian operator and the control of the corresponding eigenvalues with respect to the parameter $\r$ is more involved. In the Hermitian case let us recall the key lemma in order to control the eigenvalues with respect to a parameter:

\begin{lemma}[see \cite{EK10}]
Let $A(t)=\diag(a_1(t),\cdots,a_N(t))$ be a real {diagonal} $N\times N$-matrix and let $B(t)$ be a {\bf Hermitian} $N\times N$-matrix. Both are $C^1$ on $I\subset \R$.
Assume
\begin{itemize}
\item[(i)] $a'_j(t)\geq1$ for all $j=1,\cdots,N$ and all $t\in I$.
\item[(ii)] $||B'(t)||\leq 1/2$ for all $t\in I$.
\end{itemize}
Then
$$||(A(t)+B(t))^{-1}||\leq \frac 1 \eps $$
outside a set of $t\in I$ of Lebesgue measure {$\leq C N\eps$}.
\end{lemma}
This Lemma is false without the Hermitian hypothesis on $B$. The only way to recover a control on
$||L(\r)^{-1}||$ is to use the Cramer formula, i.e to control from above the determinant of $L(\r)$. In view of hypothesis A2 (iii) (and in particular \eqref{altern1}), we achieve this goal using the following lemma:
\begin{lemma}(see \cite{EGK, EGK1}) Let $I$ be an open interval and let $f:I \to\R$
be a $\cC^{j}$-function whose $j$:th derivative satisfies
$$\ab{f^{(j)}(x)}\ge \delta,\quad \forall x\in I.$$
Then,
$$
\Leb \{x\in I: \ab{f(x)}<\eps\}\le C (\frac\eps{\de_0})^{\frac1j}.$$
$C$ is a constant that only depends on ${j}$.
\end{lemma}
Notice that for the control of the determinant we require higher regularity with respect to the parameter $\r$ (see hypothesis A2 (iii)), the reason is the following:
if in (iii)\enskip $a\in\mathcal F$ and $b=\{\emptyset\}$, then the determinant
of $L(\r)$ is the product of $2|\F|$ term of the form $\Pi_{a\in\F}(\Om\cdot k \pm\Lambda_a+O(\eps))$ with $a\in\F$. Typically (for instance in the case of the application to the beam equation, see \cite{EGK}), we are able to prove that the first derivative of $\Om(\r)\cdot k$, in a direction depending of $k$, is large comparing to the higher derivatives. As a consequence the derivative of order $2|\F|$ of $\det L(\r)$ will be bounded from below. So in that case we take $s^*\geq2|\F|$ in hypothesis 2 (iii). For a~similar reason, if $a,b\in\mathcal F$, then we should choose $s^*\ge (2|\mathcal F|)^2$.

\subsection{Iteration}\label{iter}

In this section we would like to explain why the iteration combines a finite linear iteration
with a ``super-quadratic'' infinite iteration.

As we have seen in Section~\ref{hom} the KAM proof is based on an infinite sequence of change of variables like
$$
(h+f)\circ \Phi_S=h_++f_+,
$$
where we expect $f_+^T$ is ``small as'' $( f^T)^2$. But
actually $f_+^T$ is not really quadratic in term of $f^T$: we get (here $[\cdot]$ denotes a convenient norm)
$$[f_+^T]\sim e^{-\ga\Delta_+}[f^T]+\Delta^{\exp}e^{2\ga d_\Delta}[f^T]^2.$$
The factor $\Delta^{\exp}e^{2\ga d_\Delta}$ occurs because the diameter of the blocks $\le d_{\Delta}$ interferes
with the exponential decay and influences the equivalence between the $l^\infty$-norm and the operator-norm. The term $e^{-\ga\Delta_+}[f]$ comes from the fact that we do not solve the homological equation for blocks $[a]_{\Delta_+} \neq [b]_{\Delta_+}$ with $|a|=|b|$ (see Section~\ref{mel} and in particular \eqref{yes}).\\
So at step $k$, if $f_k^T=O(\eps_k)$, we would like
$$e^{-\ga_k\Delta_{k+1}}\eps_k+e^{2\ga_k d_{\Delta_k}} \eps_k^2\sim \eps_k^2\, .$$
This is not possible: $\ga_k d_{\Delta_k}\leq 1$ and $\ga_k\Delta_{k+1}\geq-\ln\eps_k$ are not compatible. Actually in \cite{EGK, EGK1}, as in \cite{EK10}, at each step of our infinite iteration, we apply a finite Birkhoff procedure to obtain
$[f^T_{+}] \sim [f^T]^K.$ Precisely we will choose $K_k=[\ln \eps_k^{-1}]$. The crucial fact is that, during all the $K$ Birkhoff steps, the normal form is not changed and thus the small divisors are not changed. As a consequence the clustering remains the same, i.e. $\Delta$ is fix and thus the "bad term" $e^{2\ga d_\Delta}$ is also fix. Then we iterate the previous procedure with a new clustering associated to $\Delta^+\gg \Delta$ to obtain
$$[f_+^T]\sim e^{-\ga\Delta_+}[f^T]+\Delta^{\exp}e^{2\ga d_\Delta}[f^T]^K$$
and now we will be able to control the bad terms due to the growth of the clusters.

\section{Applications}

The first part of this section is devoted to two examples that does not require a singular KAM theorem
since in both of them we use external parameters to avoid resonances. In other word, to prove the results of these examples we can use Theorem~\ref{main} with $\delta_0=1$.
 Nevertheless these two examples generalize the existing results.

 Next, in Section~\ref{sing}, we
 present the application of our main theorem to a singular situation --
 the beam equation without external parameters. We have detailed this result in \cite{EGK}, it requires the refined version, see Theorem~\ref{main2}, of our abstract KAM Theorem.

\subsection{Beam equation with a convolutive potential
}\label{s4.1}
Consider the $d$ dimensional beam equation on the torus
\be \label{beamm}u_{tt}+\Delta^2 u+V\star u + \eps g(x,u)=0 ,\quad x\in \T^{d}.
\ee
 Here $g$ is a real analytic function on $\T^{d}\times I$, where $I$ is a
 neighborhood of the origin in $\R$, and the
 convolution potential $V:\ \T^{d}\to \R$ is supposed to be analytic with
 real Fourier coefficients $\hat V(a)$, $a\in\Z^{d}$.

 Let $\cA$ be any subset of cardinality $n$ in $\Z^{d}$. We set $\L=\Z^{d}\setminus \cA$, 
$
\rho=(\hat V_a)_{a\in\cA},
$
and treat $\rho$ as a parameter of the equation,
$$
\rho=(\rho_{a_1},\dots,\rho_{a_n})
\in\D=[\rho_{a'_1},\rho_{a''_1}]\times\dots \times[\rho_{a'_n},\rho_{a''_n}]
$$
 (all other Fourier coefficients are fixed). We denote $\mu_a=|a|^4+ \hat V(a)$, $a\in\Z^{d}$,
 and assume
 that $\mu_a >0$ for all $a\in \cA$, i.e. $|a|^4+\rho_a>0$ if $a\in\cA$. We also suppose that
 $$
 \mu_l\ne0,\quad \mu_{l_1}\ne\mu_{l_2} \qquad \forall \, l, l_1, l_2 \in\L,\ \ab{l_1}\ne \ab{l_2}.
 $$

 Denote
$$
\F=\{a\in\L: \mu_a<0\}, \;\; |\F|=: {{N}},\quad \L_\infty = \L\setminus \F\,,
$$
consider the operator
$$
\Lambda=|\Delta^2 +V\star\ |^{1/2}=\diag \{\Lambda_a, a\in \Z^{d}\}\,, \quad \Lambda_a= \sqrt{|\mu_a|}\,,
$$
and
the following operator $\Lambda^{\#}$,
linear over real numbers:
$$
\LL(ze^{i \langle a, x\rangle })
= \left\{\begin{array}{ll}\
\ z\la e^{i \langle a, x\rangle }, \;\; a\in\L_{\infty}\,,\\
 -\bar z \la e^{i \langle a,x\rangle }, \;\; a\in\F,
\end{array}\right.
$$
Introducing the complex variable
 $$
 \psi= \frac 1{\sqrt 2}(\Lambda^{1/2}u- i\Lambda^{-1/2}\dot u)= (2\pi)^{-d/2} \sum_{a\in\Z^{d}}\xi_a e^{i \langle a, x\rangle }\,,
 $$
 we get for it the equation (cf. \cite[Section 1.2]{EGK})
\be\label{k1}
\dot \psi=i \big( \Lambda^{\#}\psi+ \eps \frac1{\sqrt2}\Lambda^{-1/2} g\left(x,\Lambda^{-1/2} \left(\frac{\psi+\bar\psi}{\sqrt 2}\right)\right)\,.
\ee
Writing $\xi_a=(u_a+iv_a)/\sqrt2$ we see that eq.~\eqref{k1} is a
Hamiltonian system with respect to the symplectic form $\sum dv_s\wedge du_s$ and the Hamiltonian $H=h+\eps f$, where
$$
 f= \int_{\T^{d}} G\left(x,\Lambda^{-1/2} \left(\frac{\psi+\bar\psi}{\sqrt 2}\right)\right) \dd x\,,\qquad \p_u G(x,u)=g(x,u)\,,
$$
and $h$ is the quadratic Hamiltonian
$$
h ({u}, {v}) = \sum_{a\in\cA} \Lambda_a |{\psi}_a |^2
+\Big\langle
\bH\left(\begin{array}{c} {u}_\F \\ {v}_\F \end{array}\right),
\left(\begin{array}{c} {u}_\F \\ {v}_\F \end{array}\right)
\Big\rangle
+ \sum_{a\in\L_\infty} \Lambda_a |\xi_a |^2 \,.
$$
Here
 ${u}_\F={}^t({u}_a, a\in\F)$ and
 $\bH$ is a symmetric $2{{N}}\times2{{N}}\,$-matrix. The $2{{N}}$ eigenvalues of the Hamiltonian operator with
 the matrix $\bH$ are the real numbers $\{\pm\Lambda, a\in\F\}$.
 So the linear system \eqref{beamm}${}\mid_{\eps=0}$ is stable if and only if $N=0$.

Let us fix any $n$-vector $I=\{I_a>0,a\in\cA\}$.
The $n$-dimensional torus
\ben \left\{\begin{array}{ll}
 |\xi_a|^2 =I_a,\quad &a\in \cA\\
\xi_a=0,\quad & a\in \L=\Z^{d}\setminus \cA,
\end{array}\right.
\een
is invariant for the unperturbed linear equation; it is linearly stable if and only if ${{N}}=0$.
In the linear
space span$\{\xi_a, a\in\cA\}$ we introduce the action-angle variables $(r_a,\theta_a)$ through the relations
$
\xi_a=\sqrt{(I_a+r_a)}e^{i\theta_a}$, $ a\in\cA. $
The unperturbed Hamiltonian becomes
\be \label{hbeam}
h= \text{const} +\Om(\r)\cdot r
+\Big\langle
\bH\left(\begin{array}{c} {u}_\F \\ {v}_\F \end{array}\right),
\left(\begin{array}{c} {u}_\F \\ {v}_\F \end{array}\right)
\Big\rangle
 +
 \sum_{a\in\L_\infty}\Lambda_a |\psi_a|^2\,,
\ee
with $ \Om(\r)=( \Omega_a(\r)=\Lambda_a(\r)=\sqrt{|a|^4+\r_a}, \,{a\in\cA} )$, and the perturbation becomes
\be \label{fbeam}
f=\eps \int_{\T^{d}}G
\left(x, \hat u(r,\theta;\zeta)(x)
\right)\dd x, \quad \hat u(r,\theta;\zeta)(x) = \Lambda^{-1/2}\Big(\frac{\psi+\bar\psi}{\sqrt2}\Big),
\ee
 i.e.
$$
\hat u
=\sum_{a\in\cA} \frac
{ \sqrt{(I_a+r_a)} \,(e^{i\theta_a}\phi_a + e^{-i\theta_a}\phi_{-a} )}
 {\sqrt{ 2\Lambda_a}}
 + \sum_{a\in\L}\frac{\xi_a\phi_a +\bar\xi_a\phi_{-a}}{\sqrt{ 2\Lambda_a}}.
$$
In the symplectic coordinates $((u_a, v_a), a\in \L)$ the Hamiltonian $h$ has
 the form \eqref{h},
 and we wish to apply to the Hamiltonian $h=h+\eps f$
 Theorem~\ref{main}.

 The Hypothesis~A1
 with a constant $c$ of order one and $\beta=2$ holds trivially. The
 Hypothesis~A2 also holds since for each case (i)-(iii) the second alternative
 with $\omega(\rho)=\rho$
 is fulfilled for $s^*=(2N)^2$ (see Section~\ref{hype}) and for some $\delta_0\sim1$.
 Since the discrete set $\{\Lambda_a-\Lambda_b\mid{a,b\in\L_\infty}\}$ accumulates only on the integers, Hypothesis A3 reduces essentially to a diophantine condition on $\Om(\r^*)$. As $\r\mapsto \Om(\r)$ is a local diffeomorphism at each point of $\R^n$, we verify that Hypothesis A3 holds true with $C\sim1$ and $\tau=n+1$.

 Finally, the function $f$ belongs to
 $\cT_{\ga,\vark,\D}(\s,\mu)$ with $\varkappa=1$ and suitable
 constants $\gamma_1, \gamma_2, \sigma, \mu>0$
 in view of Lemma~A.1 in \cite{EGK}. In particular the decay property on the hessian (see Remark~\ref{decay}) is a consequence of the smoothing property satisfied by the nonlinearity $f$: $
 \O^{\ga}(\s,\mu)\ni x\mapsto \nabla_\zeta f(x,\yy)\in Y_{\ga+1}$.

\medskip

 Let us set $ u_0(\theta,x) = \hat u(0,\theta;0)(x) $.
 Then for every $I\in\R_+^n$ and $\theta_0\in\T^{d}$
 the function $(t,x)\mapsto u_0(\theta_0+t\om,x)$ is a solution of \eqref{beamm} with
 $\eps=0$. Application of
 Theorem~\ref{main} gives us the following result:

\begin{theorem}\label{t72}
Fix $s>d/2$. There exist $\eps_*, \a, C>0$ such that for $0\leq\eps\leq \eps_*$ there is a Borel subset
${\D_\eps}\subset \D$,
$\, \meas(\D\setminus{\D_\eps})\leq C\eps^\a$,
such that for $\rho\in{\D_\eps}$ there is a function
$ u_1(\theta,x)$, analytic in $\theta\in\T^n_{\frac\s 2}$ and $H^{s}$-smooth in $x\in\T^{d}$, satisfying
$$\sup_{|\Im\theta|<\frac\s 2}\|u_1(\theta,\cdot)-u_0(\theta,\cdot)\|_{H^{s}(\T^{d})}
\leq C\eps,$$
and there is a mapping $\om':{\D_\eps}\to \R^n$,
$\ \|\om'-\om\|_{C^1({\D_\eps})}\leq C\eps,$
such that for $\r\in {\D_\eps}$ the function
$\
u(t,x)=u_1(\theta+t\om'(\r),x)
$
is a solution of the beam equation \eqref{beamm}.
Equation \eqref{k1}, linearised around its solution $\psi(t)$, corresponding to
the solution $u(t,x)$ above, has exactly ${{N}}$ unstable directions.
 \end{theorem}

 The last assertion of this theorem follows from the last part of Theorem~\ref{main} which
 implies that the linearised equation, in the directions corresponding to $\L$, reduces
 to a linear equation with a coefficient matrix which can be written as $B=B_\F \oplus B_\infty$.
 The operator $B_\F $ is close to the Hamiltonian operator with the matrix $H$, so it has ${{N}}$ stable
 and ${{N}}$ unstable directions, while the matrix $B_\infty$ is skew-symmetric, so it has imaginary
 spectrum.

\begin{remark} This result was proved by Geng and You \cite{GY1}
 for the case when the perturbation
 $g$ does not depend on $x$ and the unperturbed linear equation is stable.
\end{remark}

\subsection{NLS equation with a smoothing nonlinearity
}\label{nls}
Consider the NLS equation with the Hamiltonian
$$
g(u)=\tfrac12\int|\nabla u|^2\,dx+\frac{m}2\int|u(x)|^2\,dx
+\eps\int f(t,\De^{-\alpha}u(x),x)\,dx,
$$
where $m\ge0, \ \alpha>0,
\ u(x)$ is a complex function on the torus $\T^{d}$ and $f$ is a real-analytic
function on $\R\times \R^2\times\T^{d}$
 (here we regard $\C$ as $\R^2$). The corresponding
Hamiltonian equation is
\begin{equation}\label{-2.1}
\dot u=i \big(-\Delta+mu+\eps \De^{-\alpha}\nabla_2 f(t,\De^{-\alpha}u(x),x)\big)\,,
\end{equation}
where $\nabla_2$ is the gradient with respect to the second variable $u\in\R^2$.
 We have to introduce in this equation a vector-parameter $\rho\in\R^n$.
 To do this we can either assume that $f$ is time-independent and
 add a convolution-potential term $V(x,\rho)*u$ (cf. \eqref{beamm}),
 or assume that $f$ is a quasiperiodic function of time, $f=F(\rho t,u(x),x)$, where $\rho\in\D\subset\R^n$.
 Cf. \cite{Ber2}.

Let us discuss the second option. In this case the
 non-autonomous equation \eqref{-2.1} can be written as an autonomous system on the
 extended phase-space $\O\times\T^n\times L_2=\{(r,\theta,u(\cdot))\}$, where
 $ L_2=L_2(\T^{d};\R^2)$ and $\O$ is a ball in $\R^n$, with the Hamiltonian
 \begin{equation*}
 \begin{split}
&g(r,u,\rho)=h(r,u,\rho)+\eps\int F(\theta,\De^{-\alpha}u(x),x)\,dx,\\
&h(r,u, \rho)=\langle \rho, r\rangle +
\tfrac12\int|\nabla u|^2\,dx+\frac{m}2\int|u(x)|^2\,dx.
\end{split}
\end{equation*}
Assume that $m>0$ \footnote{\ if undesirable, the term $imu$ can be removed from eq.~\eqref{-2.1}
by means of the substitution $u(t,x)=u'(t,x)e^{imt}$.}
 and take for $A$ the operator $-\Delta+m$ with the eigenvalues
$\Lambda_a=|a|^2+m$. Then the Hamiltonian $g(r,u,\rho)$ has the form, required by Theorem~\ref{main}
 with
 $$
 \L=\Z^{d_*},\ \F=\emptyset, \ \varkappa=\min(2\alpha,1), \
 \beta=2,\quad c,C,\delta_0\sim 1 \text{ and } \tau=n+1 $$
 (any $\beta_2$ will do here in fact)
 and suitable $\sigma, \mu, \gamma_1>0$ and $\ga_2=m_*$ (in particular the decay property on the hessian of $f$ is a consequence of the regularization of order $2\alpha$ imposed on the nonlinearity). The
 theorem applies and implies that, for a typical $\rho$, equation \eqref{-2.1}
 has time-quasiperiodic solutions of order $\eps$. The equation, linearised
about these solutions, reduces to constant coefficients and all its Lyapunov exponents are zero.

If $\alpha=0$, equations \eqref{-2.1} become significantly more complicated. Still the
assertions above remain true since they follow from the KAM-theorem ~ in \cite{EK10}.
Cf. \cite{EK09}, where is considered nonautonomous linear Schr\"odinger equation, which is
equation \eqref{-2.1} with the perturbation $\eps (-\Delta)^{-\alpha}\nabla_2f$ replaced by $\eps V(\rho t,x)u$,
and it is proved that this equation reduces to an autonomous equation by means of a time-quasiperiodic linear
change of variable $u$. In \cite{Ber2} equation \eqref{-2.1} with $\alpha=0$ and $f=F(\rho t,\De^{-\alpha}u(x),x)$
 is considered for the case
when the constant-potential term $mu$ is replaced by $V(x)u$ with arbitrary sufficiently smooth
potential
$V(x)$. It is proved that for a typical $\rho$ the equation has small time-quasiperiodic solutions, but not that the linearised equations are reducible to constant coefficients.

\subsection{A singular perturbation problem}\label{sing}

In \cite{EGK} we apply Theorem~\ref{main} to construct small-amplitude solutions of
 the multi-dimensional beam equation on the torus:
\be\label{beam} u_{tt}+\Delta^2 u+m u = - g(x,u)\,,\quad u=u(t,x), \quad \ x\in \T^{d}.\ee
Here $g$ is a real analytic function satisfying
\be\label{g} g(x,u)=4u^3+ O(u^4).\ee
Following Section~\ref{s4.1}, the linear part becomes a Hamiltonian system with a Hamiltonian $h$ of the form \eqref{hbeam}, with $\F$ void.
$h$ satisfies Condition~A1 (for all $m>0$) and Condition~A3 (for a.a. $m>0$ with $m$-dependent
parameters $C,\tau$), but it does not satisfy Condition~A2.

The way to improve on $h$ is to use a (partial) Birkhoff normal form around $u=0$ in order to extract a piece from the non-linear part which improves on $h$. This leads to a situation where the Assumption~A2 and the size of the perturbation are linked -- a singular perturbation problem. In order to apply our KAM theorem ~to such a singular situation one needs a careful and precise description of how the smallness
requirement depends on the Assumptions~A2-A3. In other word we have to explicit \eqref{epsest2} to be able to manage the case $\delta_0\to 0$ (see Remark \ref{rem-sing}). This is quite a serious complication which is carried out in paper \cite{EGK} to obtain
\begin{theorem}\label{main2}
Assume that Hypotheses~A1-A3 hold for $\yy\in\D$. Assume that
\be\label{chi}\chi, \xi\,= O(\delta_0^{1-\aleph})\ee
for some $\aleph>0$.
Then there exist
$\eps_0, \ka, \bar\beta >0$ independent of $\delta_0$ and $\aleph$
such that if
\be\label{epsest}
\eps\big( \log \frac1\eps\big)^{\bar\beta} \le \eps_0 \delta_0^{1+\ka\aleph} =: \eps_*
\ee
then
all the statements of Theorem \ref{main} remain true.\end{theorem}
To apply this refined version of our KAM theorem, we first have to put the Hamiltonian in convenient normal form.
 Let us try to give an overview of this normal form procedure.

Let $\cA$ be a finite subset of $\Z^d$, $|\cA|=:n\ge0$. We define
$$
\L = \Z^d\setminus \cA\,.
$$
%
Let us take a vector with positive components $I=(I_a)_{a\in\cA}\in \R^n_+$.
The $n$-dimensional real torus
\ben
T^n_I=
 \left\{\begin{array}{ll}
\xi_a=\bar\eta_a,\;
 |\xi_a|^2 =I_a,\quad &a\in \cA\\
\xi_s=\eta_s=0,\quad & s\in \L \,,
\end{array}\right.
\een
is invariant for the linear hamiltonian flow (i.e. $g=0$ in \eqref{beam}).
We prove the persistency of most of the
 tori $T^n_I$
 when the perturbation $f$ (given by \eqref{fbeam}) turns on, assuming that the set of nodes $\cA$ is {\it admissible}
 or {\it strong admissible}\ in the following sense:
 For vectors $a,b \in\Z^d$ we
write
$$
a \an b \quad \text{iff} \quad \#\{x\in \Z^d \mid |x|=|a| \text{ and } |x-b| = |a-b|\} \le 2\,.
$$
Relation $a\an b$ means that
the integer sphere of radius $|b-a|$ with the centre at $b$ intersects the integer sphere $\{ x\in\Z^d\mid |x|=|a| \} $
in at most two points.
 \begin{definition}\label{adm}
A finite set $\cA\in\Z^d$, $|\cA|=:n\ge0$, is called admissible iff
$$
 a,b\in\cA, \ a\ne b
\Rightarrow |a|\neq |b|\,.
$$
An admissible set $\cA$ is called strongly admissible iff
$$
 a,b\in\cA, \ a\ne b
\Rightarrow |a|\an |b|\,.
$$
\end{definition}
Certainly if $|\cA|\leq 1$, then $\cA$ is admissible, but for $|\cA| > 1 $ this is not true. For $d \leq 2$ every admissible set is strongly admissible, but in higher dimension this is no longer true: in Appendix~B of \cite{EGK} it is proved that $\cA=\{(0,1,0),(1,-1,0)\}$ is admissible but not strongly admissible.
However, strongly admissible, and hence admissible sets are typical: see \cite{EGK} Appendix~E. \smallskip

 Now we focus on the quadratic part of the nonlinearity:
$$f_4=\int_{\T^d}u^4dx=
 (2\pi)^{-d}\sum_{(i,j,k,\ell)\in\J}\frac{(\xi_i+\eta_{-i})(\xi_j+\eta_{-j})(\xi_k+\eta_{-k})(\xi_\ell+\eta_{-\ell})}{4\sqrt{\li\lj\lk\lel}}\,,
$$
where $\J$ denotes the zero momentum set:$$\J:=\{(i,j,k,\ell)\subset\Z^d\mid i+j+k+\ell=0\}.$$
After a standard Birkhoff normal form that kills all the non resonant term, $f_4$ reduces, for generic $m$, to the resonant part
\begin{align*}
Z_4=&\frac 3 2(2\pi)^{-d}\sum_{\substack{(i,j,k,\ell)\in\J \\ \{|i|,|j| \}=\{|k|,|\ell| \}}}\frac{ \xi_i\xi_j\eta_k\eta_\ell }{\li\lj}.
\end{align*}
The terms corresponding to $(i,j,k,\ell)\in\cA^4$ will modify the internal frequencies: $\om\leadsto \Om(I)$. It is relatively simple to verify that there are no terms corresponding to exactly three indices in $\cA$ (see \cite{EGK}). It remains to consider the terms of the form
\begin{align*}P&=\xi_a\xi_b\eta_k\eta_\ell \quad \text{for }a,b\in\cA,\ \ell,k\in\L,\ (a,b,k,\ell)\in\J \text{ and } \{|a|,|b|\}=\{|k|,|\ell|\}\\
Q&= \xi_a\xi_k\eta_b\eta_\ell\quad \text{for }a,b\in\cA,\ \ell,k\in\L ,\ (a,b,k,\ell)\in\J \text{ and } \{|a|,|k|\}=\{|b|,|\ell|\} .\end{align*}
The set
$$\L_f=\{b\in\L\mid \exists a\in\cA \text{ such that } |b|=|a|\}$$
will play an important role: it corresponds to exterior modes that are resonant with some internal mode.\\
For $(k,\ell)\in \L\setminus \L_f$, there is no term of type $Q$, and the only terms of type $P$ are
 $\xi_a\eta_b\xi_k\eta_k$, they will contribute to the new frequencies $\lambda_k\leadsto\Lambda_k,\ k\in\L\setminus\L_f$. Now when $(k,\ell)\in \L_f$ the situation is more complicated and gives rise to elliptic modes and, possibly, hyperbolic
 modes: $\L_f=\L_e\cup\L_h$. Let us denote $\L_\infty=\L\setminus\L_h$.\\
In \cite{EGK} we construct an analytic
 symplectic change of variables
 $$
 \tilde\Phi_I : (r',\theta',u,v) \mapsto (r, \theta, \xi, \eta),
 $$
 such that the transformed Hamiltonian $H_I = H\circ\tilde\Phi_I$ reads
 \be\label{transff}
\begin{split}
H_I&=
\Om(I) \cdot r' +\frac12
\sum_{a\in\L_\infty} \Lambda_a(I)(u_a^2+v_a^2)\\
&+ \frac\nu2 \Bigg(\sum_{b\in\L_e} \Lambda_b(I) \Big( u_{b}^2 + v_{b}^2\Big)
 + \Big\langle \bH(I) \left(\begin{array}{ll}u^h\\ v^h \\ \end{array}\right), \left(\begin{array}{ll}u^h\\ v^h \\ \end{array}\right)
 \Big\rangle\Bigg)
+\tilde f(r',\theta', \tilde\zeta; I)
\end{split}
\ee
where $\bH(I)$ is some explicit real symmetric matrix.

It turns out that $H_I$ satisfies Hypothesis~A2 with $\delta_0\lsim |I|$ while the nonlinearity $\tilde f(\cdot;I)$ is of size $\sim|I|$ and its jet is of size $\sim |I|^{3/2}$. Therefore, taking $\delta_0 =|I|^{5/4}$ and $\aleph$ small enough (depending on $\ka$), the smallness requirement \eqref{chi}-\eqref{epsest} is satisfied for $|I|$ small enough.

Application of Theorem~\ref{main} then leads to
\begin{theorem}
There exists a zero-measure
 Borel set $\Cc\subset[1,2]$ such that for any strongly admissible set $\cA\subset\Z^d $, $|\cA|=:n\ge1$, any analytic
 nonlinearity \eqref{g} and any $m\notin\Cc$ there exists
 a Borel set $\fJ\subset \R^n_+$, having density one at the origin, with the following property:

 \noindent
 There exist
 a constant $C>0$, an exponent $\a>0$, a continuous mapping
 $\
 U: \T^n\times \fJ \to Y^R
 $
,
 analytic in the first argument, satisfying
 \be\label{dist1}
 \big| U(\T^n\times\{I\}) - ( \sqrt{I} \, e^{i \theta}\,, \sqrt{I} \, e^{-i \theta}, 0)
 \big|_{Y^R}
 \le C |I|^{1 -\a_*}\
 \ee
 and a continuous
 vector-function
 \be\label{dist11}
 \om' : \fJ \to \R^n\,,
\ee
 such that
 for any $I\in \fJ$ and $\theta\in\T^n$ the parametrized curve
\be\label{solution}
 t \mapsto U(\theta +t\omega'(I), I)
\ee
is a solution of the beam equation \eqref{beam}.
 Accordingly, for each $I\in \fJ$ the analytic $n$-torus $U(\T^n\times\{I\}) $
 is invariant for equation~\eqref{beam}. \\
 Furthermore this invariant torus is linearly stable if and only if $\L_h=\emptyset$.
\end{theorem}

The torus $T^n_I$ is invariant for the linear beam equation \eqref{beam}${}_{g=0}$. For $m\notin\Cc$ and $I\in\fJ$
 the constructed invariant torus $U(\T^n\times\{I\})$
 of the nonlinear beam equation is a small perturbation of $T^n_I$.
 \smallskip

Denote $\cT_\cA=U(\T^n\times\fJ)$. This set is invariant for the beam equation and is filled in with its time-quasiperiodic solutions.
Its Hausdorff dimension equals $2|A|$. Now consider $\cT = \cup \cT_\cA$, where the union is taken over all strongly admissible sets $\cA \subset\Z^d$. This invariant set has infinite Hausdorff dimension. Some time-quasiperiodic solutions of (1.1), lying on $\cT$, are linearly stable, while, if $d \geq 2$, then some others are unstable.
For instance, it is proved in \cite{EGK}, Appendix B, that the invariant
 2-tori, constructed on the strongly admissible
set $\cA=\{(0,1),(1,-1)\}$ are unstable.

\end{document}